# Quantization Opportunities for Polyharmonic Subdivision Wavelets Applied to Astronomical Images

Ognyan Kounchev, Damyan Kalaglarsky

*Abstract:* We continue the study of a new family of multivariate wavelets which are obtained by "polyharmonic subdivision". We provide the results of experiments considering the distribution of the wavelet coefficients for the Lena image and for astronomical images. The main purpose of this investigation is to find a clue for proper quantization algorithms.

## 1 Introduction

The big success of JPEG relies upon an ingenious quantization algorithm, [9]. Let us recall that in the standard JPEG algorithm which is based on $8\times 8$ rectangle pixel tiling of an image, we normally assume that the image has a size which is the power of $2$. For example, we take $128\times 128$ or $256\times 256$. To every rectangle $R_j$ for $j = 1,2,...,N$, one applies the Discrete Cosine Transform (DCT). Hence, one obtains $64$ coefficients $\{R_k^j\}_{k=1}^{64}$ which are put in one array, and proceeds with a **quantization algorithm**. An important point is the statistical observation about the distribution of these $N\times 64$ coefficients; in the case of a $128\times 128$ pixel image, $N = 256$, hence we have a total of $256\times 64 = 16384$ coefficients.

Alternatively, one may apply a wavelet transform to the $8\times 8$ pixel or bigger tiles of the image. This is the idea of JPEG2000 where the tiling of the image is not into $8\times 8$ pixels but adaptively chosen, and a wavelet transform is applied to the rectangles $R_j$. However the performance of JPEG2000 is reported by different authors to be no more than $20\%$ better than JPEG. To be more precise, let the one-dimensional signal $f$ have the expansion

$$f(t) = \sum_{k\geq K} f_{k,j}\psi_{k,j}(t) + \sum g_{k,j}\varphi_{k,j}(t) = D + A$$

where we have taken several levels of details *D*; here $\psi_{k,j}$ are the mother wavelets, and their coefficients $\{f_{k,j}\}$ are starting with the biggest "detail" until level $K$, and respectively $\varphi_{k,j}$ are the dilations-shift of the father wavelet which are representing the Approximation *A* below level $K$. Roughly speaking, the terms containing the mother wavelets $f_{k,j}\psi_{k,j}(t)$ capture the "edge structures" of the image, while the terms containign the father wavelets $g_{k,j}\varphi_{k,j}(t)$ capture the "intra-edge structures" of the image. Both sets of coefficients $\{f_{k,j}\}$ and $\{g_{k,j}\}$ enjoy a nice "wavelet type" structure for natural images: very few of them are really large while the majority are small. However one needs to find more structure in the smaller coefficients which would allow for better quantization and compression of the natural images. This is especially urgent in the analysis of images. Indeed, the efforts in this direction by a generation of people devoted to Wavelet Analysis show, that up to now, it does not seem that the existing wavelet families possess such a structure susceptible to quantization. Beyond the standard tensor product wavelets one has to mention also a more tricky constructions having elongated support which are able to capture "directional information", as curvelets and ridgelets (Donoho-Candes), bandlets (Mallat), contourlets (Vetterli), brushlets (Coifman), shearlets (Kutyniok), and others, cf. e.g. [6], [7].

In the present paper we continue the study of the new family of wavelets which have elongated support and were termed Polyharmonic Subdivision Wavelets, [3], [4], [5].

The **main purpose** of the present work is to investigate experimentally the possibility to find structures in the coefficient sets of the Polyharmonic Subdivision Wavelets which yield to quantization.

## 2 Definitions and the wavelet algorithm

Let $N \geq 1$ be a fixed integer. It will correspond to the integer $N$ in the construction of Daubechies (cf. [1], formula (6.6)). We have the following definitions and the algorithm:

1. For fixed $\xi \geq 0$ and $k \in \mathbb{Z}$ we define the polynomial

$$d(z) := d^{[k],\xi}(z) := d^{[k]}(z) := \frac{(z+x_0)^N (z^{-1}+x_0)^N}{(1+x_0)^{2N}} \quad \text{for } z \in \mathbb{C}; \tag{1}$$

here we put $x_0 := e^{-\xi 2^{k+1}}$ and

$$\eta = \eta^{[k],\xi} := \frac{4x_0}{(1+x_0)^2} = \frac{2}{1+\cosh(\xi/2^{k+1})}.$$

2. We define the polynomial

$$R_N(x) := \sum_{j=0}^{N-1} \binom{N+j-1}{j} y^j \tag{2}$$

(which corresponds to the polynomial $P_N$ in the notations of [2] and [1]), and define the key polynomial

$$Q(x) = (2-\eta)^{-N} R_N\left(\frac{1-\eta(1-x)}{2-\eta}\right). \tag{3}$$

3. Further, we define the trigonometric polynomial $b^{[k]}(z)$ by putting

$$b^{[k]}(z) := b^{[k],\xi}(e^{i\omega}) := Q^{[k],\xi}\left(\sin^2\frac{\omega}{2}\right) \tag{4}$$

where we have the notations

$$x = \sin^2\frac{\omega}{2} = \frac{1-\cos\omega}{2} = \frac{1}{2} - \frac{z+z^{-1}}{4}.$$

4. Let us denote the zeros of the polynomial (2) by $c_j^D$, i.e.

$$R_N(y) = \sum_{j=0}^{N-1} \binom{N+j-1}{j} y^j = \frac{(2N-2)!}{((N-1)!)^2} \prod_{j=1}^{N-1}(y-c_j^D).$$

Then the polynomial $Q$ as defined in formula (3) is given by

$$Q(x) = (2-\eta)^{-2N+1} \eta^{N-1} \frac{(2N-2)!}{((N-1)!)^2} \prod_{j=1}^{N-1}(x-C_j),$$

where

$$C_j := \frac{c_j^D(2-\eta)+\eta-1}{\eta}. \tag{5}$$

5. We take the trigonometric polynomial

$$M_1(z) := \frac{(z+x_0)^N}{(1+x_0)^N} \tag{6}$$

as the " square root" of $d^{[k]}(z)$, i.e. $d^{[k]}(z) = |M_1(z)|^2$ for $|z|=1$.

6. The polynomial $M_2$ of degree $\leq N-1$ such that

$$|M_2(e^{i\omega})|^2 = \frac{1}{2}Q\left(\sin^2\frac{\omega}{2}\right), \tag{7}$$

is obtained by using the roots $c_j^D$ of the polynomial $R_N$. Let the polynomial $Q$ have the zeros $C_j$ defined in (5) and let us put $c_j = 1 - 2C_j$. We see that $Q\left(\sin^2\frac{\omega}{2}\right) = \vec{Q}(\cos\omega)$ for some polynomial $\vec{Q}$ and $c_j$ are the zeros of $\vec{Q}$. At this point we apply the algorithm for the Riesz representation of $\vec{Q}$, see e.g. [1], p. 197–198.

7. For every integer $m \geq 0$ and a real number $\xi \geq 0$ the family of functions

$$M(z) := M^{[m]}(z) := M^{[m],\xi}(z) := M_1(z)M_2(z) \tag{8}$$

represents the refinement masks for the family of scaling functions (father wavelets) $\{\varphi_m(t)\}_{m\geq 0}$. Let us remark that they are different at different levels since the Multiresolution Analysis is non-stationary.

8. The construction of the mother wavelets $\{\psi_m(t)\}_{m\geq 0}$ (and their filter coefficients) does not differ from the standard construction for the usual stationary wavelets, as in [2], [1], [8]. In particular, for the scaling function $\varphi_m$ we have

$$\varphi_m(t) = \sum_{j\in Z} M_j^{[m],\xi} \varphi_m(2t-j),$$

Hence, the mother wavelet is obtained as

$$\psi_m(t) = \sum_{j\in Z} (-1)^j M_{1-j}^{[m],\xi} \varphi_m(2t-j).$$

Here $M_j^{[m],\xi}$ are the coefficients of the polynomial $M(z)$.

9. For simplicity we assume that the Image which we consider is a function $u(t,y)$, for $t \in [a,b]$ and is $2\pi$-periodic in $y \in R$. We have

$$u(t.y) = \sum_{\eta \in Z^n} v_\eta(t) e^{i\langle \eta, y\rangle} \quad \text{where the Fourier coefficients are} \tag{9}$$

$$v_\eta(t) = \frac{1}{2\pi}\int_0^{2\pi} u(t.y) e^{-i\langle \eta, y\rangle} dy.$$

For every fixed $\eta \in Z$ and $\xi = |\eta|$ we make Wavelet Analysis which is based on the non-stationary father and mother wavelets $\{\varphi_m^\xi(t)\}_{m\geq 0}$ and $\{\psi_m^\xi(t)\}_{m\in Z}$.

10. Assume that we fix some level of the "Details" of the wavelet expansion. Thus we have the following expansion

$$u(t.y) = \sum_{m\geq m_0, j\in Z}\left(\sum_{\eta\in Z}\gamma_{\eta,m,j}\psi_m^{|\eta|}(t-j)e^{i\langle\eta,y\rangle}\right) + \sum_{m<m_0, j\in Z}\left(\sum_{\eta\in Z}g_{\eta,m,j}\varphi_m^{|\eta|}(t-j)e^{i\langle\eta,y\rangle}\right) \tag{10}$$

$$= D + A \tag{11}$$

This means that we perform the wavelet expansion only down to a certain level $m_0$ in the details part $D$. The "Approximation" part $A$ contains the expansion in father wavelets.

## 3  Experiments
### 3.1  Experiments with the Lena image

In Figure 1 below we provide the distributions of the (complex) coefficients $\{\gamma_{\eta,m,j}\}$ and $\{g_{\eta,m,j}\}$ for the Lena image of $128\times128$ pixels. On the other hand, due to boundary effects the coefficients are not precisely the same number $128^2$ but their number is about $26.000$ in the complex plane:

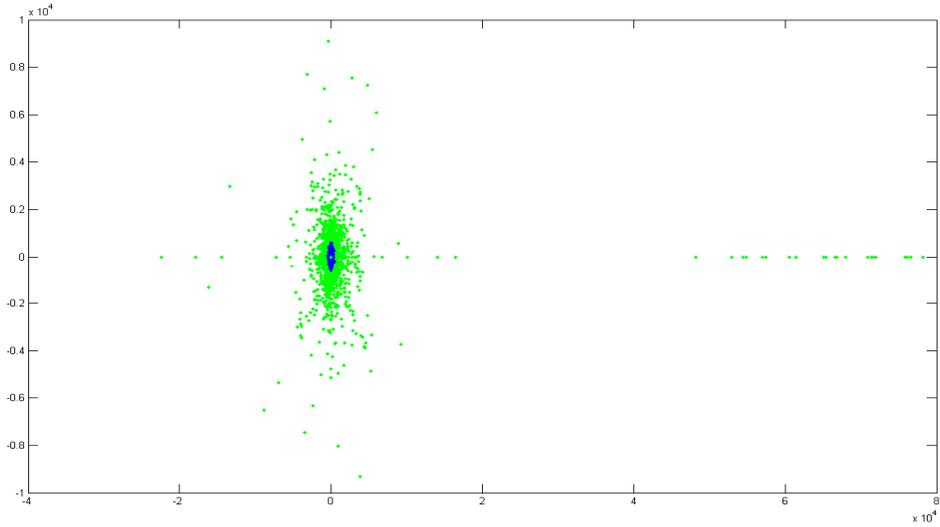

Figure 1

We see that there are few coefficients which are really large of size about $20.000$. The coefficients have a very intricate structure which is revealed only after zooming them. The following two figures, Figure 2 and Figure 3, show this structure:

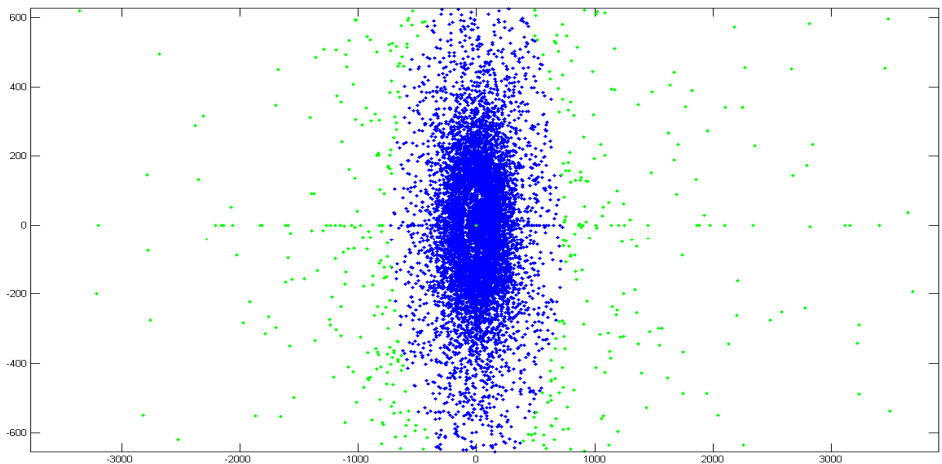

Figure 2

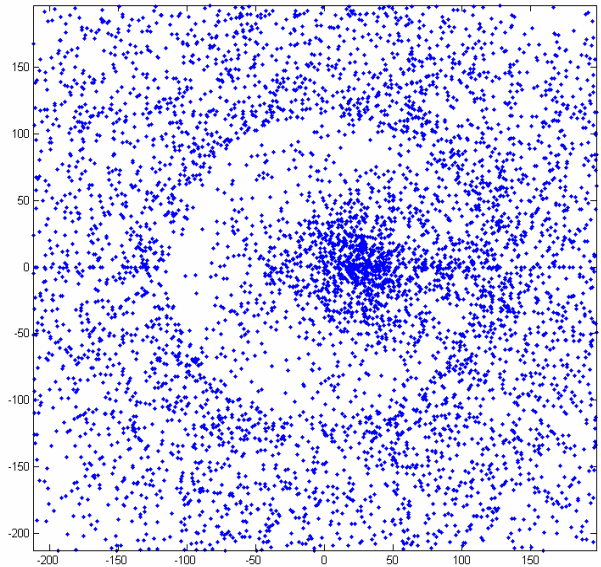

Figure 3

On the other hand, it is instructive to see what is the distribution of the coefficients of the usual two-dimensional Daubechies wavelets on Figure 4.

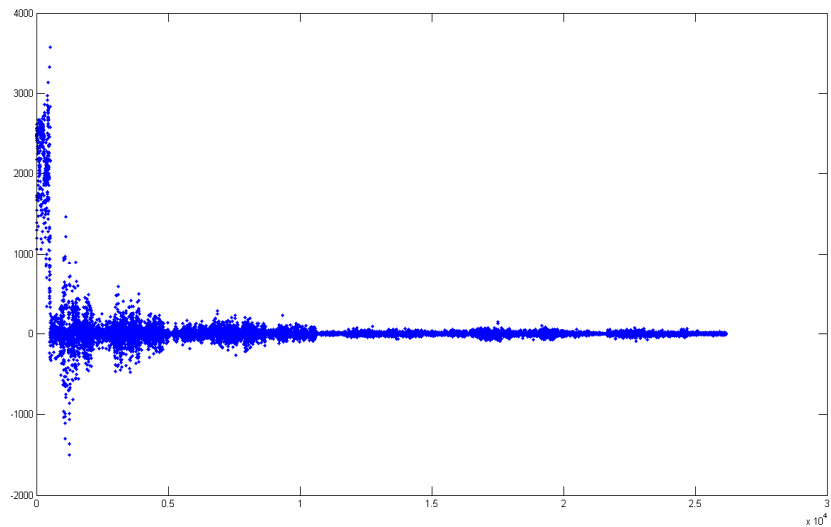

Figure 4

Finally, we have the approximation of the Lena image in Figure 5 below. On the top right we have the PH wavelets and at the bottom are DB (the usual Daubechies wavelets). We have the following parameters of the application of the methods:
   *Original image*: Lena128_Greyscale.bmp
   *Reconstructed images*: (upper left - original, upper right - PH, lower left - DB);
   Used: 9bit quantizer
   **Compression ratio of PH Wavelet** encoding without rotation $= 3.1196$
   **PSNR of PH Wavelet** reconstruction after quantization $= 80.7606$

**Compression ratio of DB Wavelet** encoding = 3.0910
**PSNR of DB Wavelet** reconstruction after quantization = 72.7275

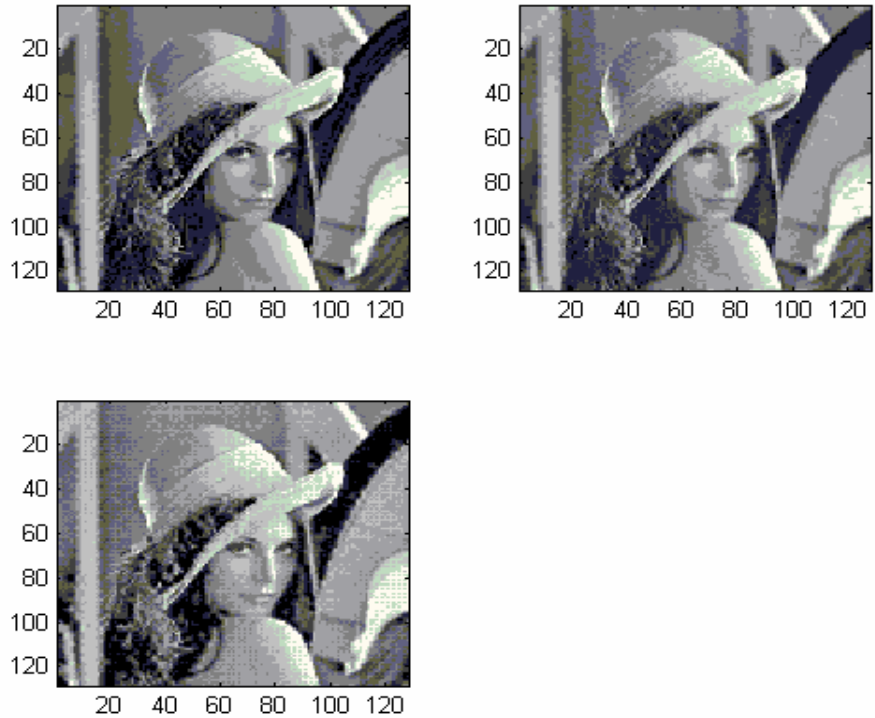

Figure 5

### 3.2 Experiments with Astronomical images

We use the image **ROZ050 000046** of Rozhen Observatory, which is represented by a $128 \times 128$ pixel image.

In Figure 6 we obtain the distribution of the coefficients $\{\gamma_{\eta,m,j}\}$ and $\{g_{\eta,m,j}\}$:

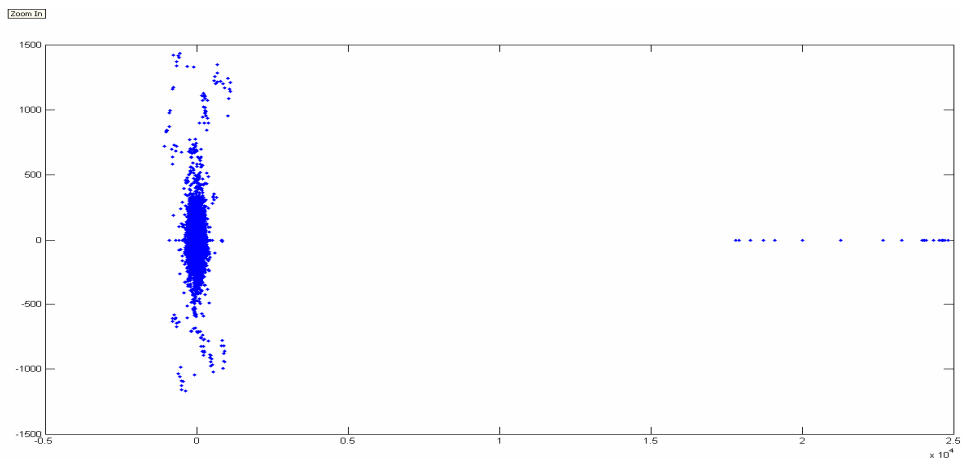

Figure 6

As in the case of the Lena image we see that the coefficients are very much concentrated in a small area. It is important to emphasize that these coefficients are more concentrated than the corresponding coefficients in the Lena image. It is interesting however that the coefficient set has a structure similar to the one we have discovered in the case of the Lena image, and again it may be discovered only after zooming it; we provide the following two Figures 7 and 8 where we see this structure:

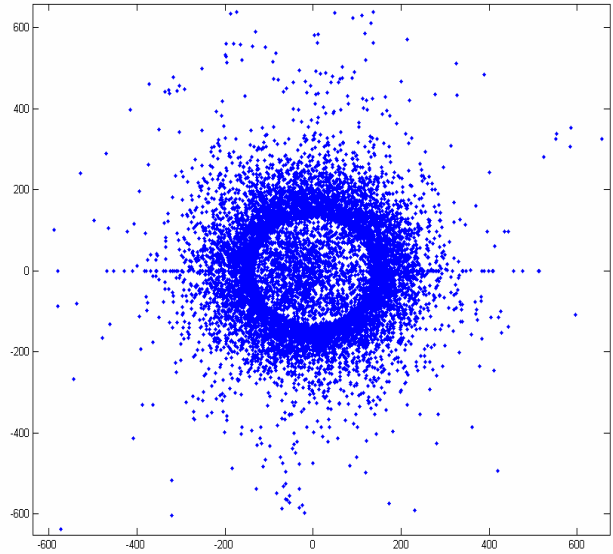

Figure 7

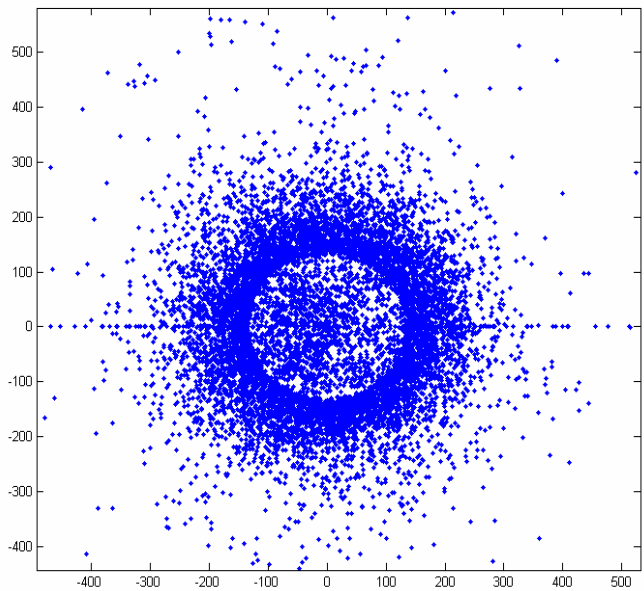

Figure 8

In Figure 9 we see completely clearly that there is a circle pattern close to the origin:

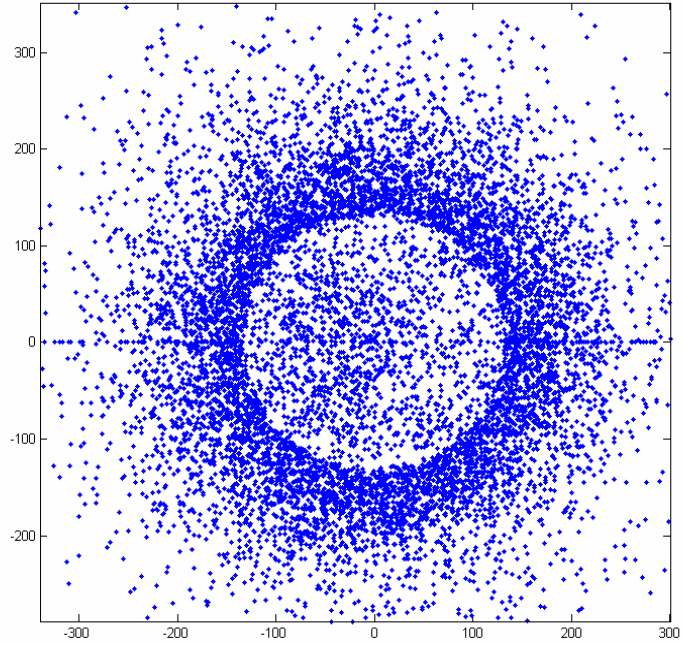

Figure 9

In Figure 10 we provide the coefficients of the application of the DB wavelets to the Astronomical image. We see that they are more scattered to a large area compared to the coefficients in the case of the Lena image:

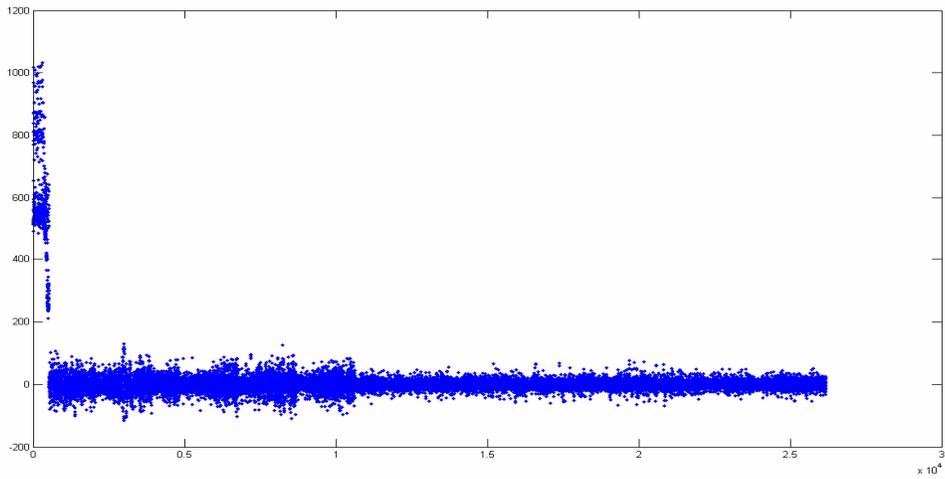

Figure 10

In Figures 11 and 12 we have the original image of the astronomical plate with the result of application of the PH wavelets and the standard DB wavelets:

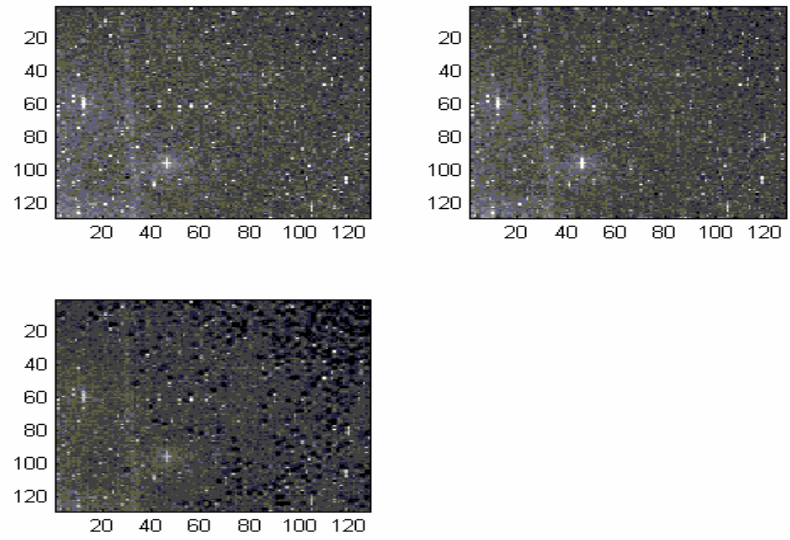

Figure 11

and the following which has less compression ratio:

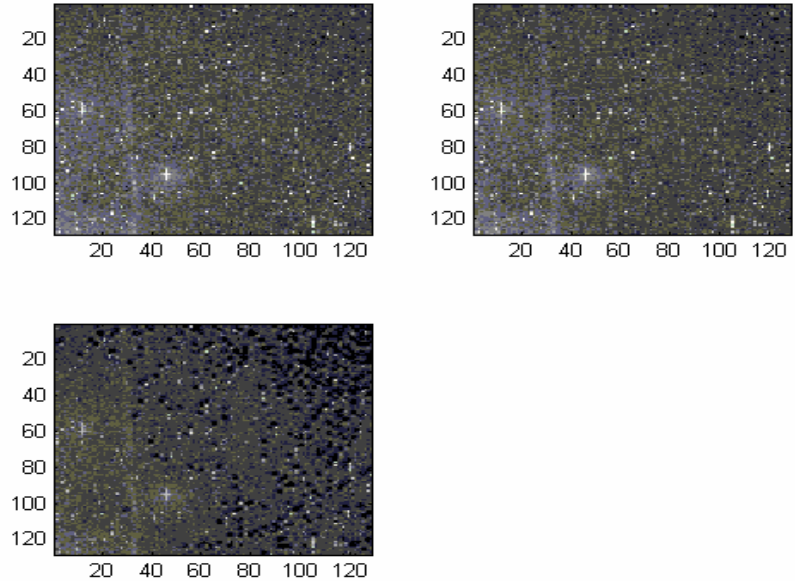

Figure 12

The parameters of the application are the following:
Original image: ROZ050 000046_128x128_vert.gif
For Figure 11:
**Compression ratio of PH Wavelet** encoding without rotation = 2.7012

**PSNR of PH Wavelet** reconstruction after quantization $= 81.1902$
**Compression ratio of DB Wavelet** encoding $= 2.4414$
**PSNR of DB Wavelet reconstruction** after quantization $= 78.1313$
The parameters for Figure 12 above are:
**Compression ratio of PH Wavelet** encoding without rotation $= 2.0208$
**PSNR of PH Wavelet reconstruct** after quantization $= 83.2997$
**Compression ratio of DB Wavelet** encoding $= 2.0231$
**PSNR of DB Wavelet reconstruction** after quantization $= 78.4700$

**Conclusion** *The coefficients of the Polyharmonic Subdivision Wavelets show non-trivial structure in the complex plane which may be susceptible to successful application of further quantization procedures.*

**ABOUT THE AUTHORS**
Prof. Ognyan Kounchev, PhD, D.Sc., Institute of Mathematics and Informatics, Bulgarian Academy of Science & IZKS, University of Bonn; Phone: +359 2 979 3851, E-mail: kounchev@gmx.de.
Damyan Kalaglarsky, Institute of Astronomy, Bulgarian Academy of Science, Phone: +359 2 979 3851, E-mail: damyan@skyarchive.org.